\documentclass[a4paper,10pt]{amsproc}
\usepackage{amsfonts,lingmacros,tree-dvips, amsmath, amssymb, graphicx, mathrsfs}

\usepackage[all]{xy}

\newdir{ >}{!/6pt/@{ }*:(1,0.2)@_{>}*:(1,-0.2)@^{>}}

\theoremstyle{plain}

\newtheorem{theorem}{Theorem}[section]

\newtheorem{proposition}[theorem]{Proposition}
\theoremstyle{definition}
\newtheorem{definition}[theorem]{Definition}

\newtheorem{counter example}[theorem]{Counter Example}

\numberwithin{equation}{section}

\DeclareMathAlphabet{\mathscr}{OT1}{pzc}{m}{it} 
\begin{document}
\Large{
		\title{ON NON-BAIRE SETS IN CATEGORY BASES}
		
		\author[S. Basu]{Sanjib Basu}
		\address{\large{Department of Mathematics,Bethune College,181 Bidhan Sarani}}
		\email{\large{sanjibbasu08@gmail.com}}
		
		\author[A. Deb Ray]{Atasi Deb Ray}
		\address{\large{Department of Pure Mathematics, University of Calcutta, 35, Ballygunge Circular Road, Kolkata 700019, West Bengal, India}}
		\email{\large{debrayatasi@gmail.com}}
		
		\author[A.C.Pramanik]{Abhit Chandra Pramanik}
		\address{\large{Department of Pure Mathematics, University of Calcutta, 35, Ballygunge Circular Road, Kolkata 700019, West Bengal, India}}
		\email{\large{abhit.pramanik@gmail.com}} 
		
		\thanks{The third author thanks the CSIR, New Delhi – 110001, India, for financial support}
	\begin{abstract}
	In this paper, we prove a result on non-Baire sets in category bases which when applied together with a result of Grzegorek yeilds a comparatively stronger version of a decomposition theorem due to Ulam.
	\end{abstract}
\thanks{}
	\maketitle

\section{INTRODUCTION}
In $[2]$ the following two theorems were proved
\begin{theorem}
	There corresponds to every non-measurable subset M of the real line, a subclass $\langle$E$\rangle$ of the class of all measurable sets of positive Lebesgue measure (having the same cardinality as that of the later) such that M$\cap$E is non-measurable for every E $\in$ $\langle$E$\rangle$.
\end{theorem}
\begin{theorem}
 There corresponds to every non-Baire subset B of the real line, a subclass $\langle$E$\rangle$ of the class of all second category sets having Baire property(having the same cardinality as that of the later) such that B$\cap$E is non-Baire for every E$\in$$\langle$E$\rangle$.
\end{theorem}
The real line with its usual topology is second countable and this fact plays an effective role in proving the above two theorems. But topological spaces, let alone metric spaces in general are not second countable, so in order to extend the above two results in higher settings we need to develop some appropriate generalizations of second countability. Now second countability can be generalized by the aid of `measure zero cardinal' in measure theoretic situations or by the use of $\sigma$-locally finite base in topological situations and these were used in $[1]$ to give extensions of the above two theorems. But no extension which could unify the above two theorems was given in $[1]$ and is probably not possible without using the notion of category bases. Here in this paper, we first give an unification of Theorem $1.1$ and Theorem $1.2$ in the setting of category bases. We then apply this result and the Grzegorek's unification $[3],[9]$ of two theorems of Sierpi\'{n}ski $[4],[5],[9]$ to prove a comparatively stronger version of a decomposition theorem due to Ulam $[9]$.
\section{PRELIMINARIES AND RESULTS}
The idea of category base which is a generalization of both measure and topology and whose main objective is to present measure and category(topology) and also some other aspects of point set classification within a common framework was introduced by J.C.Morgan in the mid seventies of the last century and has developed since then through a series of papers $[6],[7],[8]$ etc.
We recall some of the basic definitions and results in this subject which may be found in the above references and also in $[9]$.
\begin{definition}
A category base is a pair (X,$\mathcal{C}$) where X is a set and $\mathcal{C}$ is a family of subsets of X, called regions satisfying the following set of axioms:
\begin{enumerate}
\item Every point of X belongs to some region; i,e., X=$\cup$$\mathcal{C}$
\item Let A be a region and $\mathcal{D}$ be a non-empty family of regions having cardinality less than the cardinality of $\mathcal{C}$\\
i) If A$\cap$($\cup$$\mathcal{D}$) contains a region, then there is a region D$\in$$\mathcal{D}$ such that A$\cap$D contains a region.\\
ii) If A$\cap$($\cup$$\mathcal{D}$) contains no region, then there is a region B$\subseteq$A that is disjoint from every region in $\mathcal{D}$.
\end{enumerate}
\end{definition}
\begin{definition}
In a category base (X,$\mathcal{C}$), a set is called singular if every region contains a subregion which is disjoint from the set. Any set which can be expressed as countable union of singular sets is called meager. Otherwise, it is called abundant.
\end{definition}
\begin{definition}
	In a category base (X,$\mathcal{C}$), a set S is called Baire if in every region, there is a subregion in which either S or its complement X$-$S is meager.
\end{definition}
Consequently it follows that a set is non-Baire if there is a region in every subregion of which both the set and its complement are abundant.
\begin{proposition}
	If (X,$\mathcal{C}$) is a category base, $\mathcal{N}$ is a subfamily of $\mathcal{C}$ with the property that each region in $\mathcal{C}$ contains a region in $\mathcal{N}$ and Y=$\cup$$\mathcal{N}$, then (Y,$\mathcal{N}$) is also a category base and the $\mathcal{N}$-singular sets coincide with $\mathcal{C}$-singular subsets of Y. In addition, if U is a subset of Y and Y$-$U is $\mathcal{N}$-singular, then X$-$U is $\mathcal{C}$-singular.
\end{proposition}
\begin{proposition}
	If (Y,$\mathcal{N}$) is a category base, then there exists a disjoint subfamily $\mathcal{M}$ of $\mathcal{N}$ such that Y$-$$\cup$$\mathcal{M}$ a singular set. Moreover, $\mathcal{M}$ may be so selected that for every region N, there exists a region M$\in$$\mathcal{M}$ such that N$\cap$M contains a region.
\end{proposition}
\begin{theorem}
	Every abundant set in a category base (X,$\mathcal{C}$) is abundant in every subregion of some region.\\
	The above theorem is known as the ``Fundamental Theorem".
\end{theorem}
We now formulate an unification of Theorem $1.1$ and Theorem $1.2$ .To start with, we introduce
\begin{definition}
	A set A in a category base (X,$\mathcal{C}$) has the property ($\star$) if there exists a region C such that for every subregion D of C and every set E which essentially contains D , A$\cap$E is non-Baire.\\
	(A set E contains D essentially means that D$-$E is meager in the category base (X,$\mathcal{C}$)) 
\end{definition}
\begin{theorem}
	In any category base, every non-Baire set satisfies the ($\star$)-property.
	\begin{proof}
	Let us assume that there is a non-Baire set A which does not have the ($\star$)-property. This means that for every C$\in$$\mathcal{C}$, there exists D$\in$$\mathcal{C}$ such that D$\subseteq$C and a set E such that D is essentially contained in E and A$\cap$E is Baire. We may assume that E$\subseteq$D.
	
	The collection of all such regions D taken together constitutes a subcollection $\mathcal{N}$ of $\mathcal{C}$ such that every region in $\mathcal{C}$ contains a region in $\mathcal{N}$. Let Y=$\cup$$\mathcal{N}$ and by Proposition 2.5 we can extract out a subfamily $\mathcal{M}$ (of $\mathcal{N}$) of mutually disjoint regions having the property stated therein. Consider enumerations \{$D_\alpha$: $\alpha$$<$$\Theta$\} (of $\mathcal{M}$), \{$E_\alpha$: $\alpha$$<$$\Theta$\} of the above collections where $\Theta$ is the smallest ordinal whose cardinality is same as the cardinality of $\mathcal{M}$$-$\{$\phi$\}. Let D*=$\bigcup\limits_{\alpha<\Theta}$$D_\alpha$, E*=$\bigcup\limits_{\alpha<\Theta}$$E_\alpha$. Clearly, E*$\subseteq$D*. By construction, Y$-$D* is $\mathcal{N}$-singular and so by Proposition 2.4, X$-$D* is $\mathcal{C}$-singular.
	
	The set A$\cap$E* is Baire. Take any region D$\in$$\mathcal{N}$. Then there exists a region $D_\alpha$$\in$$\mathcal{M}$ such that D$\cap$$D_\alpha$ contains a region $G_\alpha\in\mathcal{N}$. Since A$\cap$$E_\alpha$ is Baire, $G_\alpha$ contains a subregion $H_\alpha\in\mathcal{C}$ such that either A$\cap$$E_\alpha$ or its complement in $H_\alpha$ is meager. Since A$\cap$E*=$\bigcup\limits_{\alpha<\Theta}$ (A$\cap$$E_\alpha$) where the sets A$\cap$$E_\alpha$ are mutually disjoint, this implies that either A$\cap$E* or its complement in $H_\alpha$ is meager.But every region in  $\mathcal{C}$ contains a region in  $\mathcal{N}$. So A$\cap$E* is Baire.
	
	The set D*$-$E* is meager. Since $D_\alpha$$-$$E_\alpha$ is meager, we may write  {$D_\alpha$}$-${$E_\alpha$}=$\bigcup\limits_{n=1}^{\infty}P_{n}^{(\alpha)}$ ($\alpha$$<$$\Theta$) where each $P_{n}^{(\alpha)}$ is singular. We choose n arbitrarily and fix it and using a procedure similar as above, find that every region contains a subregion {$H_\alpha$}$\in$$\mathcal{C}$  which is disjoint from  $P_{n}^{(\alpha)}$ and hence from $\bigcup\limits_{n=1}^{\infty}P_{n}^{(\alpha)}$. As every region in $\mathcal{C}$ contains a subregion in $\mathcal{N}$, the set $P_{n}$=$\bigcup\limits_{\alpha<\Theta}P_{n}^{(\alpha)}$ is $\mathcal{C}$-singular. Hence D*$-$E*=$\bigcup\limits_{\alpha<\Theta}(D_{\alpha}-E_{\alpha})$	=$\bigcup\limits_{\alpha<\Theta}{\bigcup\limits_{n=1}^{\infty}P_{n}^{(\alpha)}}$ =$\bigcup\limits_{n=1}^{\infty}{\bigcup\limits_{\alpha<\Theta}P_{n}^{(\alpha)}}$ =$\bigcup\limits_{n=1}^{\infty}P_{n}$ is  $\mathcal{C}$-singular. Therefore A is a Baire set in (X, $\mathcal{C}$) which is a contradiction.
	\end{proof}
\end{theorem}
\begin{definition}
	A category base (X, $\mathcal{C}$) is called point-meager if every singleton set in it is meager. It is a Baire base if every region is abundant.
\end{definition}
The following two theorems were proved by Sierpi\'{n}ski $[4], [5], [9]$. Here $\aleph_{0}$ and $\aleph_{1}$ represent the first infinite and the first uncountable cardinals.
\begin{theorem}
	If $2^{\aleph_{0}}$=$\aleph_{1}$, then every subset of $\mathbb{R}$ which is of second category in every interval contains an uncountable family of disjoint sets each of which is of second category in every interval.
\end{theorem}
\begin{theorem}
	If $2^{\aleph_{0}}=\aleph_{1}$ ,then every subset of $\mathbb{R}$ with positive Lebesgue outer measure contains an uncountable family of disjoint sets each of which has the same measure as the given set.
\end{theorem}
Grzegorek [3], [9] unified the above two theorems in point-meager,Baire bases.
\begin{theorem}
	Let (X,$\mathcal{C}$) be a point-meager,Baire base. If $\mathcal{C}$ satisfies CCC (countable chain condition) and every region has cardinality $\aleph_{1}$, then every abundant set S can be decomposed into an uncountable family of disjoint sets each of which is abundant in every region in which S is abundant.
\end{theorem}
We now show that using Theorem 2.8 and Theorem 2.12, we can derive under some modification of hypothesis, a comparatively stronger version of the following decomposition theorem of Ulam.  
\begin{theorem}
	Let (X,$\mathcal{C}$) be a point-meager base. Then every set of cardinality $\aleph_{1}$ which is not a Baire set can be decomposed into an uncountable family of sets none of which is a Baire set.
\end{theorem}

To every non-Baire set there is associated a region that determines its non-Baireness in the sense that in every subregion of that region both the set and its complement are abundant. We call this region a test region corresponding to the non-Baire set.\\

Now from Theorem 2.12, it follows that  every non-Baire set S in (X,$\mathcal{C}$) can be decomposed into a family \{$S_\alpha\}_{\alpha < \Omega}$ ($\Omega$ is the first ordinal of cardinality $\aleph_{1}$) consisting of mutually disjoint sets $S_\alpha$ each of which is abundant in every region in which S is abundant. But S also satisfies property ($\star$) by Theorem 2.8. Hence there exists a region in every subregion of which every $S_\alpha$ and its complement are abundant proving that 
\begin{theorem}
	If (X,$\mathcal{C}$) is a point-meager, Baire base where $\mathcal{C}$ satisfies CCC and every region has cardinality $\aleph_{1}$, then every non-Baire set of cardinality $\aleph_{1}$ can be decomposed into an uncountable family of non-Baire sets which is an uniformly non-Baire family in the sense that there is a common test region for all the members in the family.
\end{theorem}
This property of uniform non-Baireness of the family does not follow from the proof of Ulam's theorem in any point-meager base. Now since every abundant set contains a non-Baire set, so from the above theorem it also follows that 
\begin{theorem}
	If (X,$\mathcal{C}$) is a point-meager, Baire base where $\mathcal{C}$ satisfies CCC and every region has cardinality $\aleph_{1}$, then every abundant set of cardinality $\aleph_{1}$ can be decomposed into an uncountable family of non-Baire sets which is an uniformly non-Baire family in the sense that there is a common test region for all the members in the family .
\end{theorem}
Now depending on the following definition and result, Theorem 2.14 and Theorem 2.15 may be expressed in another equivalent form.
\begin{definition}
	In a category base (X,$\mathcal{C}$), a set A is called completely non-Baire in a region D if for every Baire set B such that B$\cap$D is abundant, both A$\cap$B and (D$-$A)$\cap$B are abundant.
\end{definition}
\begin{theorem}
	In any Baire base, a set is non-Baire if and only if it is completely non-Baire in some region.
	\begin{proof}
		Let A be a non-Baire set in a Baire base (X,$\mathcal{C}$). Then according to the definition, there is a region D in $\mathcal{C}$ in every subregion of which both A and its complement (X$-$A) are abundant. Let B be any Baire set such that B$\cap$D is abundant. Since every region is here a Baire set by hypothesis, we may assume that B$\subseteq$D. Again B being Baire, from the Fundamental theorem, it follows that there is a subregion C of D which is essentially contained in B. Consequently, both B$\cap$A and B$\cap$(D$-$A) are abundant.
		
		Conversely, suppose A is completely non-Baire in a region D. Since by hypothesis every region is an abundant Baire set, both A and its complement are abundant in every subregion of D. Hence A is non-Baire.
	\end{proof}
\end{theorem}
\begin{theorem}
	If (X,$\mathcal{C}$) is a point-meager, Baire base where $\mathcal{C}$ satisfies CCC and every region has cardinality $\aleph_{1}$, then every non-Baire set of cardinality $\aleph_{1}$ can be decomposed into an uncountable family of non-Baire set each of which is completely non-Baire in a region which can be chosen independently of the choice of members in the family.
\end{theorem}
\begin{theorem}
	If (X,$\mathcal{C}$) is a point-meager,Baire base where $\mathcal{C}$ satisfies CCC and every region has cardinality $\aleph_{1}$, then every abundant set of cardinality $\aleph_{1}$ can be decomposed into an uncountable family of non-Baire sets each of which is completely non-Baire in a region which can be chosen independently of the choice of the members in the family.
\end{theorem}
\textbf{NOTE:} From Theorem 2.8 it follows that in any Baire base (X,$\mathcal{C}$) there corresponds to every non-Baire set A, a subclass  $\langle$E$\rangle$ of the class of abundant Baire sets such that A$\cap$E is non-Baire for every E$\in$ $\langle$E$\rangle$. An explicit description of this subclass is also given. But in this general situation, it is not possible to give an indication as to what the cardinality of this subclass should be. However, if the category base is a perfect base[9], then since there exists in it a singular set having the cardinality of continuum (Th 36,Ch 5,[9]), the cardinality of the subclass will be same as the cardinality of the class of all abundant Baire sets.

\bibliographystyle{plain}

	\end{document}